\documentclass{degm}
\usepackage{amsmath,math,amsthm,epsfig}
\newcommand\circuits{{\mathcal C}}
\newcommand\prim{{\mathcal P}}
\newcommand\rels{{\mathcal R}}
\newcommand\scount{{\operatorname{sc}}}
\newcommand\SC{{\mathcal O}}
\newcommand\proba{{\mathbb P}}

\begin{document}
\title{Cactus Trees, and Lower Bounds on the Spectral Radius of
  Vertex-Transitive Graphs}
\author{Laurent Bartholdi}
\date{July 2nd, 2002}
\markboth{Laurent Bartholdi --- May 20, 2002}{Cactus Trees, and
  Estimations of the Spectral Radius of Vertex-Transitive Graphs}
\maketitle
\begin{abstract}
  This paper gives lower bounds on the spectral radius of
  vertex-transitive graphs, based on the number of ``prime
  cycles'' at a vertex. The bounds are obtained by constructing
  circuits in the graph that resemble ``cactus trees'', and
  enumerating them. Counting these circuits gives a coefficient-wise
  underestimation of the Green function of the graph, and hence an
  underestimation of its spectral radius.
  
  The bounds obtained are very good for the Cayley graph of surface
  groups of genus $g\ge2$, with standard generators (these graphs are
  the $1$-skeletons of tessellations of hyperbolic plane by $4g$-gons,
  $4g$ per vertex). We have for example for $g=2$
  \[0.662418\le\|M\|\le 0.662816,\]
  and for $g=3$
  \[0.552773\le\|M\|\le 0.552792.\]
\end{abstract}

\section{Introduction: Groups}
Throughout this paper, $\Gamma$ will be a group generated by a finite,
symmetric set $S$, of cardinality $\#S=d$. ``Symmetric'' means that
$S=S^{-1}$. Many of the objects we define will depend heavily on the
choice of $S$, even though we will not make it explicit in the
notation.

The \emph{Cayley graph} $\gf$ of $\Gamma$ is the graph with vertex set
$\Gamma$, and vertices connected under the right action of $S$; i.e.\ 
$\gamma$ and $\gamma s$ are joined for all $\gamma\in\Gamma$ and $s\in
S$. The \emph{Markov operator} $M:\ell^2(\Gamma)\to\ell^2(\Gamma)$ is
defined by
\[(Mf)(\gamma)=\frac1d\sum_{s\in S}f(\gamma s).\]
It is used to study the simple random walk on $\Gamma$; for instance,
the probability of return in $n$ steps is $p_n=\langle
M^n\delta_1|\delta_1\rangle$, where $\delta_1$ is the Dirac function
at $1\in\Gamma$.

The first step in understanding $M$ is the computation of its
(operator) norm $\|M\|$, also called the \emph{spectral radius} of
$\gf$; indeed, the probabilities $p_n$ satisfy
$\limsup_{n\to\infty}\sqrt[n]{p_n}=\|M\|$. Harry Kesten showed
in~\cite{kesten:rwalks,kesten:amen} the group-theoretical importance
of $\|M\|$: we always have
\begin{equation}\label{eq:kesten}
  \frac{2\sqrt{d-1}}d\le\|M\|\le1,
\end{equation}
with equality on the left if and only if $\gf$ is a tree, i.e.\ 
$\Gamma$ is a free product of $\Z/2$'s and $\Z$'s, with $S$ consisting
of the standard generators and their inverses; and equality holds on
the right if and only if $\Gamma$ is amenable.

Assume now that $\Gamma$ is not free; say it has a relation of length
$k\ge3$. William Paschke obtained in~\cite{paschke:norm} the estimate
\begin{equation}\label{eq:paschke}
  \|M_\Gamma\|\ge\min_{s>0}\left\{2\cosh(s)+(d-2)Q\left(\frac{\cosh(ks)+1}{\sinh(s)\sinh(ks)}\right)\right\},
\end{equation}
where $Q(t)=(\sqrt{t^2+1}-1)/t$.

The purpose of this paper is to show that the lower
bound~\eqref{eq:paschke} on $\|M\|$ can be improved if additional
hypotheses are made on the number of relations of $\Gamma$, and on the
number of distinct cyclic permutations of these relations.
\begin{defn}
  The \emph{Green function} of $\Gamma$ is the formal power series
  \[G(t) = \sum_{w\in S^*:\,w\equiv_\Gamma1}t^{|w|}=\sum_{n\ge0}p_nd^nt^n,\]
  i.e.\ the growth series of the words representing $1$ in $\Gamma$.
  
  For two power series $G(t)=\sum_{n\ge0}g_nt^n$ and
  $H(t)=\sum_{n\ge0}h_nt^n$, define $G\precsim H$ to mean $g_n\le h_n$
  for all $n\ge0$.
  
  A \emph{prime relator} is a word $w\in S^*$ such that
  $w\equiv_\Gamma1$, and such that $v\not\equiv_\Gamma1$ for all
  proper subwords $v$ of $w$.
  
  A set $\rels$ of prime relators satisfies the \emph{small
    cancellation condition} $\SC(\eta)$ if for any $w,w'\in\rels$,
  and any factorization $w=uv$ and $w'=vu'$ we have either $u=u'$ or
  $|v|\le\eta\cdot\min\{|w|,|w'|\}$.
\end{defn}
Note that a prime relator is necessarily a cyclically freely reduced
word. The symbol ``$\SC$'' stands for ``overlap''. It is a notion
close, but strictly weaker than the $C(\eta)$ in small cancellation
theory.

The main result of this paper, stated for finitely generated groups,
is the following. See in Subsection~\ref{subs:main} the more general
form stated for vertex-transitive graphs:
\begin{cor}\label{cor:main}
  Let $\Gamma$ be a group generated by a finite symmetric set $S$ of
  cardinality $d$, and let $M$ be its Markov operator. Assume that
  $\Gamma$ has a set $\rels$ of prime relators satisfying a small
  cancellation condition $\SC(\eta)$.
  
  Let $f(t)=\sum_{w\in\rels}t^{|w|}$ be the growth series of $\rels$,
  and let $\zeta$ satisfy $\zeta-1+(1-1/d)\zeta
  f(\zeta^{\eta-1}/(d-1))=0$. Construct the power series
  \begin{align*}\tag{*}\label{eq:main}
    h(t)&=\frac{2(d-1)}{d-2+d\sqrt{1-4(d-1)t^2}},\\
    H(t,u)&=\frac{1-(1-u)^2t^2}{1+(1-u)(d-1+u)t^2}h\left(\frac {t}{1+(1-u)(d-1+u)t^2}\right),\\
    g_1(t,u)&=H(\zeta t,u),\\
    g_2(t)&=g_1\left(t\frac{d-f(t)}{d-(d-1)f(t)},\frac{(d-2)f(t)}{d-f(t)}\right),\\
    g_3(t)&=h(t)g_2\left(\frac{1-\sqrt{1-4(d-1)t^2}}{2(d-1)t}\right).
  \end{align*}
  Then the Green function of $\gf$ satisfies $G\succsim g_3$.
  
  Let $\rho$ be the radius of convergence of $g_3$. Then
  $\|M\|\ge1/(d\rho)$.
\end{cor}

Note that Corollary~\ref{cor:main} does not supersede Paschke's
result~\eqref{eq:paschke}, in that the bound it gives for the group
$(\Z/k)*(\Z/2)*\dots*(\Z/2)$ is inferior to Paschke's. It does,
however, give a superior bound for many groups, and in particular
surface groups.

\subsection{Surface groups}\label{subs:surface}
Consider the fundamental group of a surface of genus $g\ge2$
\[\Gamma_g = \big\langle a_1,b_1,\dots,a_g,b_g\big|[a_1,b_1]\dots[a_g,b_g]=1\big\rangle.\]
This group is non-elementary hyperbolic, hence non-amenable; its
spectral radius is $O(g^{-1/2})$ for large $g$. Simple lower and upper
bounds come respectively from $\Gamma_g$ being a quotient of a free
group of rank $2g$, and containing by Magnus' Freiheitssatz a
$(2g-1)$-generated free subgroup $\langle
a_1,b_1,\dots,a_{g-1},b_{g-1},a_g\rangle$:
\[\frac{\sqrt{4g-1}}{2g}\le\|M_g\|\le\frac{\sqrt{4g-3}+2}{2g};\]
For more details, see~\cite{b-c-c-h:growth}. Many improvements on
these bounds were obtained,
see~\cite{paschke:norm,cherix-v:1rel,b-c-c-h:growth,zuk:norm,nagnibeda:upperbd,bartholdi-c:growth}.
Currently, the the best known bounds are
\begin{thm}
  The spectral radius of the surface groups of genus $2$ and $3$ satisfy
  \begin{align*}
    0.662418\le&\|M_2\|\le 0.662816,\\
    0.552773\le&\|M_3\|\le 0.552792.
  \end{align*}
\end{thm}
These upper bounds are due to Tatiana
Nagnibeda~\cite{nagnibeda:upperbd}.


\subsection{Reduction to graphs}
Equations~\eqref{eq:kesten} and~\eqref{eq:paschke} hold more generally
for any vertex-transitive graph, i.e.\ for any graph whose group of
automorphisms acts transitively on vertices. Remember
from~\cite{paschke:norm} that there are vertex-transitive graphs with
no simply vertex-transitive automorphism group --- for instance, the
$1$-skeleton of a dodecahedron. Such a graph cannot be the Cayley
graph of a group.

In the context of $d$-regular vertex-transitive graphs, $p_nd^n$ is
the number of closed paths of length $n$ at any fixed vertex in $\gf$,
and $d\|M\|$ is the asymptotic exponential growth rate of these
numbers of closed paths.

In just the same way as any group is a quotient of a free group, any
$d$-regular graph is covered by the $d$-regular tree. Since closed
paths in the tree remain closed in the quotient graph, the spectral
radius of any $d$-regular graph is bounded from below by the spectral
radius of the $d$-regular tree; this proves the left inequality
of~\eqref{eq:kesten}.

Similarly, any vertex-transitive, $d$-regular graph with a loop of
length $k$ at each vertex is covered by the graph $P_{k,d}$ obtained
from a $k(d-2)$-regular tree by replacing each vertex by a $k$-gon and
equidistributing the edges on the $k$-gon's
vertices~\cite[Proposition~2.4]{paschke:norm}. The spectral radius of
$P_{k,d}$ can be computed using
\def\0{\cite[Theorem~9.2]{bartholdi:cogrowth}}
\begin{thm}[\cite{woess:free}; \cite{cartwright-s:free}; \0]
  Let $\mathcal X_1,\mathcal X_2$ be vertex-transitive graphs, and let
  $\mathcal X$ be their free product. Let $G_1(t),G_2(t)$ and
  $G(t)$ be their corresponding Green functions. Then
  \[\frac1{(tG(t))^{-1}}=\frac1{(tG_1(t))^{-1}}+\frac1{(tG_2(t))^{-1}}-\frac1t,\]
  where $F^{-1}(t)$ denotes the formal inverse, i.e.\ the series
  $E(t)$ such that $E(F(t))=F(E(t))=t$.
\end{thm}
Indeed if we take $\mathcal X_1$ a $d-2$-regular tree and
$\mathcal X_d$ a $k$-cycle, then $\mathcal X$ is $P_{k,d}$.

Clearly if $F(t)$ is an algebraic series, then $F^{-1}(t)$ is
algebraic of same degree, because $P(t,F(t))=0$ implies
$P(F^{-1}(t),t)=0$; also, $\deg(F_1(t)+F_2(t))\le\deg F_1(t)\deg
F_2(t)$. Simple computations show that for $\mathcal X_1$ and
$\mathcal X_2$ as above $G_1$ is algebraic of degree $2$ and $G_2$ is
algebraic of degree $(k+1)/2$. It follows that $\|M_P\|$ is an
algebraic number of degree at most $k+1$.

Another computations shows that $\|M_P\|$ is the right-hand side
of~\eqref{eq:paschke}; the inequality then follows.

Note that Paschke's estimate is valid for all vertex-transitive
graphs; the inequality is obtained by constructing a cover for all
graphs containing a $k$-cycle. We note \emph{a posteriori} that
$P_{k,d}$ is the Cayley graph of $(\Z/k)*(\Z/2)*\dots*(\Z/2)$, with
$d-2$ copies of $\Z/2$, but there is no \emph{a priori} reason for the
graph of smallest norm, $P_{k,d}$, to be the Cayley graph of a group.
Our method can also be understood as constructing a transitive graph
of minimal norm satisfying some conditions on its cycles; however,
this graph of minimal norm will not be the Cayley graph of a group.

\section{Preliminaries: Graphs}
A \emph{graph} $\gf$ is a pair $V,E$ of sets called respectively
\emph{vertices} and \emph{edges}, with maps $\alpha,\omega:E\to V$
called respectively \emph{start} and \emph{end}, and an involution
$\overline\cdot:E\to E$ such that $\alpha(\overline e)=\omega(e)$. The
\emph{degree} of a vertex is $\deg(v)=\#\{e\in E:\,\alpha(e)=v\}$. The
graph is \emph{$d$-regular} if each vertex has degree $d$.

A \emph{path} is a sequence $p=(p_1,\dots,p_n)$ of edges, with
$\omega(p_i)=\alpha(p_{i+1})$ for all $i\in\{1,\dots,n-1\}$. Its
\emph{length} is $|p|=n$, and its \emph{start} and \emph{end} are
$\alpha(p)=\alpha(p_1)$ and $\omega(p)=\omega(p_n)$. It is a
\emph{circuit} if $\alpha(p)=\omega(p)$. Paths are multiplied by
concatenation; therefore in this definition a graph is nothing but a
small $*$-category with object set $V$ and arrow set $E$.

A \emph{spike} in a path $p$ is an index $i\in\{1,\dots,|p|-1\}$ such
that $p_i=\overline{p_{i+1}}$. A path is \emph{reduced} if it has no
spike. The \emph{spike count} of a path $p$ is the number $\scount(p)$
of spikes in $p$. A circuit is \emph{prime} if it is reduced,
non-trivial, and $\alpha(p_i)\neq\alpha(p_j)$ for all $i\neq j$.

The \emph{Markov operator} on $\gf$ is the operator
$M:\ell^2(V,\deg)\to\ell^2(V,\deg)$ given by
\[(Mf)(v)=\frac1{\deg(v)}\sum_{\{e\in E:\,\alpha(e)=v\}}f(\omega(e)).\]
It governs the behaviour of the simple random walk on $\gf$: given two
vertices $v,w\in V$ the probability of going from $v$ to $w$ in $n$
steps is $\langle M^n\delta_w|\delta_v\rangle$, where $\delta_x$ is
the Dirac function at $x\in V$.

The graph $\gf$ is \emph{transitive} if for any two vertices $v,w\in
V$ there is a graph automorphism mapping $v$ to $w$. Transitivity
implies that the probability of return $p_n=\langle
M^n\delta_v|\delta_v\rangle$ does not depend on $v$; its exponential
rate of decay is $\|M\|$, the spectral radius of the random walk.

If $\Gamma$ is a group with symmetric generating set $S$, its
\emph{Cayley graph} $\gf$ has vertex set $\Gamma$ and edge set
$\Gamma\times S$, with $\alpha(\gamma,s)=\gamma$ and
$\omega(\gamma,s)=\gamma s$ and $\overline{(\gamma,s)}=(\gamma
s,s^{-1})$. Clearly, $\gf$ is a transitive graph under the left action
of $\Gamma$.

Let $\gf$ be a transitive graph, and let $\prim$ be a set of prime
circuits in $\gf$ at a given vertex $v$. By transitivity, there exists
a translate $\prim_w$ of $\prim$ at any vertex $w$. A \emph{spread} of
$\prim$ is the union of the $\prim_w$ for all $w\in V$.

The set $\prim$ satisfies the \emph{small cancellation condition}
$\SC(\eta)$ if it has a spread $\widehat\prim$ such that for any
$p,p'\in\widehat\prim$, and any factorization $p=qr$ and $p'=rq'$ we
have either $q=q'$ or $|r|\le\eta\cdot\min\{|q|,|q'|\}$.

The \emph{growth series} of a set $\prim$ of paths is the formal power
series $f(t)=\sum_{p\in\prim}t^{|p|}$. The \emph{Green function} of
$\gf$ is the growth series of the set of circuits at an irrelevant but
fixed vertex.

\subsection{Main result}\label{subs:main} The definitions given above
were tailor-cut to make Corollary~\ref{cor:main} a direct consequence
of the following result on graphs:
\begin{thm}\label{thm:main}
  Let $\gf$ be a vertex-transitive, $d$-regular graph, and let $\prim$
  be a set of prime circuits at a fixed vertex in $\gf$ satisfying
  $\SC(\eta)$. Let $f(t)$ be the growth series of $\prim$.  Construct
  the power series $g_1(t),\dots g_3(t)$ as in~\eqref{eq:main},
  Corollary~\ref{cor:main}.  Then the Green function of $\gf$
  satisfies $G\succsim g_3$.
  
  Let $\rho$ be the radius of convergence of $g_3$. Then
  $\|M\|\ge1/(d\rho)$.
\end{thm}

The main idea behind Theorem~\ref{thm:main} is the construction of
\emph{cactus trees} in $\gf$. A \emph{cactus tree} is a circuit in
$\gf$ that is built in three stages: at the first stage, a ``trunk''
is constructed, i.e.\ a circuit that freely reduces to the trivial
circuit. This trunk may not contain any subword of length
$(1-\eta)|p|$ of any prime circuit $p$.

At the second stage, ``fruits'', i.e.\ prime circuits, are inserted at
all vertices of the circuit, in such a way that the resulting circuit
is reduced.

At the third stage, ``spikes'', i.e.\ circuits that freely reduce to
the trivial circuit, are inserted at all vertices of the cactus tree.

Here is an example of cactus, with the three stages of construction
indicated in solid bold, dashed medium and solid thin lines:
\[\epsfig{file=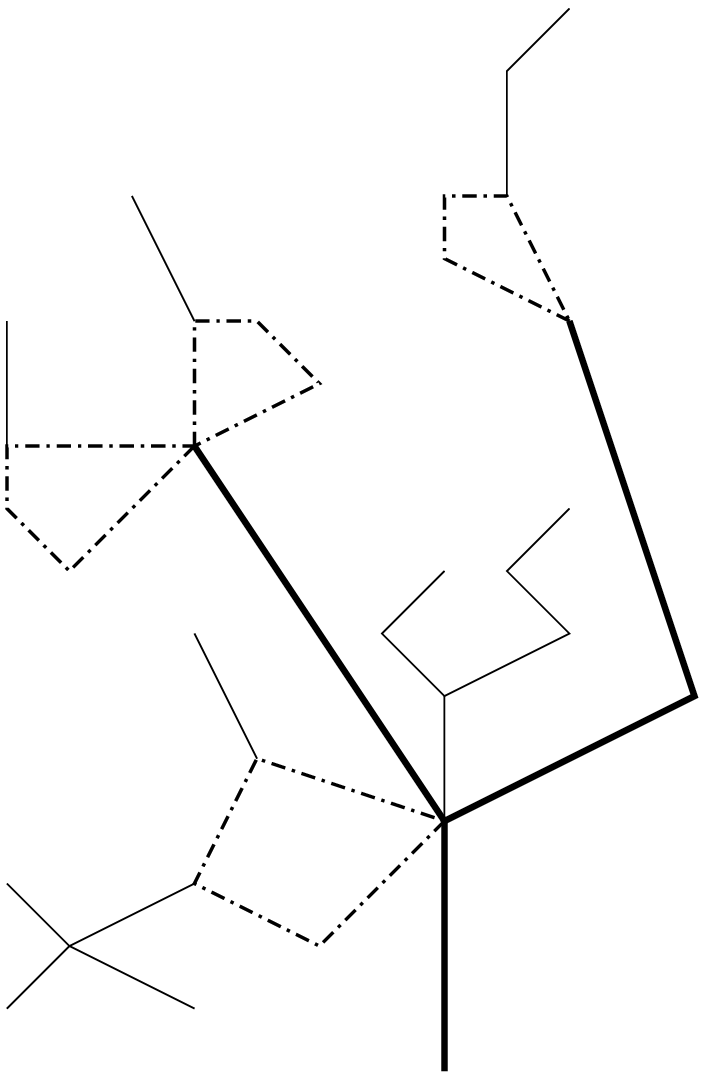}\]
Consider the graph $X_{d,m}$ that is the $1$-skeleton of the
tessellation of the hyperbolic plane by $m$-gons, $d$ per vertex.
Numerically, the growth of cactus trees is a very good underestimation
of the growth of circuits in the graph $X_{d,m}$. The Cayley graph of
the surface group $\Gamma_g$ is $X_{4g,4g}$. This explains the bounds
given in Subsection~\ref{subs:surface}.

\section{A Formula: Cogrowth}
This section recalls the main result of~\cite{bartholdi:cogrowth}. Fix
a vertex $v\in V$, and consider a set $\circuits$ of reduced circuits
at $v$. The \emph{saturation} $\langle\circuits\rangle$ of $\circuits$
is the closure of $\circuits$ under the iterated operation of
inserting spikes in paths; i.e., $\langle\circuits\rangle$ is the
smallest set of circuits containing $\circuits$ and such that for all
product of paths $pq\in\langle\circuits\rangle$ and for all $e\in E$ with
$\alpha(e)=\omega(p)=\alpha(q)$ we have $pe\overline
eq\in\langle\circuits\rangle$. (We allow $p,q$ to be empty, in which
case by convention $\omega(p)=\alpha(q)=v$).

Let $\mathcal P$ be a set of paths. Its \emph{spiky growth series} is
the formal $2$-variable power series
\[G_{\mathcal P}(t,u) = \sum_{p\in\mathcal P}t^{|p|}u^{\scount(p)}.\]
The specialization $G_{\mathcal P}(t,1)$ ``forgets'' the number of
spikes in the paths while ``remembering'' only their lengths, and in
case $\mathcal P$ consists of all the circuits at $v$ is often called
the \emph{Green function} of $\gf$ at $v$. The specialization
$G_{\mathcal P}(t,0)$ counts only reduced paths in $\mathcal P$.

\begin{thm}[\cite{bartholdi:cogrowth}, Corollary~2.6]\label{thm:cogrowth}
  Let $\circuits$ be an arbitrary set of reduced paths in a
  $d$-regular graph $\gf$, and let $\langle\circuits\rangle$ be its
  saturation. Then
  \[\frac{G_{\langle\circuits\rangle}(t,u)}{1-(1-u)^2t^2} =
  \frac{G_{\langle\circuits\rangle}\left(\frac{t}{1+(1-u)(d-1+u)t^2},1\right)}{1+(1-u)(d-1+u)t^2}.\]
  In particular, we have
  \[\frac{G_\circuits(t,0)}{1-t^2} =
  \frac{G_{\langle\circuits\rangle}\left(\frac{t}{1+(d-1)t^2},1\right)}{1+(d-1)t^2}.\]
\end{thm}

Assume $\gf$ is transitive and $d$-regular, and recall $p_n$ denotes
the probability of return in $n$ steps of the simple random walk on
$\gf$. Then at each vertex there are $d^np_n$ circuits of length $n$,
and we have the following connection between spectral radius and
counting of paths:
\begin{prop}\label{prop:lowerbd}
  Let $\circuits$ be a set of circuits at $v\in V$ in a $d$-regular
  graph $\gf$, and let $\rho$ be the the convergence radius of
  $G_{\langle\circuits\rangle}(t,1)$. Then $\|M\|\ge 1/(d\rho)$.
\end{prop}

\subsection{Forbidden words}\label{subs:fw}
Let $F$ be a set of words over an alphabet $A$, and consider the
problem of estimating the number of words over $A$ not containing any
element of $F$ as a subword. We have the following lemma, which is in
essence the Lovasz local lemma:
\begin{lem}\label{lem:local}
  Let $W$ be either $A^*$ or the set of reduced words $A^*\setminus
  \{A^*aa^{-1}A^*\}_{a\in A}$, if $A$ has an involution
  $a\leftrightarrow a^{-1}$.
  
  Let $F\subset W$ be a set of ``forbidden'' words, and set
  $L=W\setminus WFW$ be the language of words in $W$ not containing an
  element of $F$ as a subword. Let $\Phi$ denote the growth series of
  $F$, and let $\lambda,\rho$ denote the growth rate of $L,W$
  respectively; $\rho$ is either $\#A$ or $\#A-1$.
  
  Then we have $\lambda\ge\rho\zeta$, where $\zeta$ satisfies the
  equation $\zeta-1+\frac{\rho\zeta}{\#A}\Phi(\rho^{-1}\zeta^{-1})=0$.
\end{lem}
\begin{proof}
  Put on the $S_n=W\cap A^n$ the uniform distribution, and define
  events $Q_j,R_j$ on $S_n$ as follows: $Q_j$ is the set of words
  $w\in S_n$ for which there are no factorizations $w=efg$ with
  $|ef|=j$ and $f\in F$; in other words, $Q_j$ is the event of not
  containing a forbidden word ending at index $j$. Set $R_j=Q_1\cap
  Q_2\cap\dots\cap Q_{j-1}$.  Note
  $\proba(R_{n+1})=\proba(Q_n|R_n)\proba(R_n)$.
  
  We claim that $\proba(Q_j|R_j)\ge\zeta$ for all $j\in\{1,\dots,n\}$,
  and proceed by induction, the basis being
  $\proba(Q_0|R_0)=1\ge\zeta$ by definition of $\zeta$:
  \begin{align*}
    \proba(Q_j|R_j)&\ge1-\sum_{f\in F}\frac{\proba(\text{``occurence of $f$ ending at
        index $j$''}\cap R_j)}{\proba(R_j)}\\
    &\ge1-\sum_{f\in F}\frac{\proba(\text{``occurence of $f$ ending at index
      $j$''}\cap R_{j-|f|+1})}{\proba(Q_{j-1}|R_{j-1})\dots\proba(Q_{j-|f|+1}|R_{j-|f|+1})\proba(R_{j-|f|+1})}\\ 
    &\ge1-\sum_{f\in
      F}\frac{\#A^{-1}\rho^{1-|f|}\quad\proba(R_{j-|f|+1})}{\zeta^{|f|-1}\proba(R_{j-|f|+1})}\\
    &=1-\rho\zeta\Phi(\rho^{-1}\zeta^{-1})/\#A=\zeta.
  \end{align*}
  (In the third equality, we use $\proba(\text{``occurence of $f$
    ending at index $j$''})=\#A^{-1}\rho^{1-|f|}$ and independence
  with $R_{j-|f|+1}$.)

  Therefore $\#(L\cap S_n)=\rho^n\proba(R_{n+1})\ge(\rho\zeta)^n$.
%
%
%
\end{proof}

As an application, consider a $d$-regular transitive graph $\gf$, and
a spread $\prim$ of ``forbidden paths'' in $\gf$. Fix a vertex $*$ in
$\gf$, and let $T$ denote the set of circuits at $*$. By
Theorem~\ref{thm:cogrowth}, the spiky growth series of $T$ is
\[\Theta(t,1-u)=\frac{2(d-1)(1-u^2t^2)}{(d-2)(1+u(d-u)t^2)+d\sqrt{(1+u(d-u)t^2)^2-4(d-1)t^2}}.\]
\begin{lem}\label{lem:localtree}
  If $F$ is a set of freely reduced words, then with the notation of
  Lemma~\ref{lem:local} the spiky growth series of $T\cap L$ is
  coefficient-wise at least $\Theta(\zeta t,u)$.
\end{lem}
\begin{proof}
  For $n\in\N$ consider $T\cap A^n$; then $\zeta^n$ is a lower bound
  of the probability, for a word of length $n$ chosen with uniform
  probability, to belong to $L$. For any $i<j\in\{1,\dots,n\}$, if
  there exists $w\in T\cap A^n$ whose subword $w_i\dots w_j$ is
  reduced, then the restriction $T\cap A^n\to A^{j-i+1}$ to indices
  $i\dots j$ yields a uniform reduced word of length $j-i+1$. If
  follows that
  \[\zeta^n\le\mathbb P(w\in F|w\in A^n)\le\mathbb P(w\in F|w\in T\cap A^n).\]
\end{proof}

Note that in case a more detailed description of $F$ is available,
strengthenings of the above result are possible. With surface groups
in mind, we consider the concrete case
$A=\{a_1^{\pm1},b_1^{\pm1},\dots,a_g^{\pm1},b_g^{\pm1}\}$; $W$ either
$A^*$ or the set $F_A$ of reduced words over $A$; and
$r=[a_1,b_1]\dots[a_g,b_g]$.  We take for $F$ the set of subwords of
the cyclic permutations of $r^{\pm1}$ of length $4g-1$.  As above we
write $\rho$ for the growth rate of $W$. Let $\Lambda(t)$ denote the
growth series of $L$. Then
\begin{equation}\label{eq:localsg}
  \Lambda\succsim\frac1{1-\rho t+8gt^{4g-1}\frac{1-t}{1-t^{4g-1}}}.
\end{equation}
Indeed, Lemma~\ref{lem:local} applies, and can be slightly improved by
letting $F_i$ denote the set of subwords of length $i$ of $rrr\dots$
and $r^{-1}r^{-1}r^{-1}\dots$, for all $r\in F$, and writing $\Phi_i$
for the generating series of $F_i$. The equation defining $\Lambda$ is
then
\[\Lambda=\frac1{1-\rho t}\Big(1-\big(\Phi_{4g-1}-\Phi_{4g}+\Phi_{8g-2}-\Phi_{8g-1}+-\dots\big)\Lambda\Big);\]
indeed a word in $W$ is either in $L$, or is of the form $fgh$, with
$h\in L$ and $g\in F_i$ for $i$ maximal. The result follows by
inclusion-exclusion.

\section{Main Result: Proof of Theorem~\ref{thm:main}}
Assume throughout this section that a vertex $v\in V$ has been fixed
in the transitive, $d$-regular graph $\gf$. By
Proposition~\ref{prop:lowerbd}, a lower bound on $\|M\|$ can be
obtained by evaluating $G_\circuits(t,1)$ for some set $\circuits$ of
circuits at $v$.

Let $\prim$ be a set of prime circuits at $v$ satisfying the small
condition condition $\SC(\eta)$, and let $f=G_\prim(t,1)$ be its
growth series.  Let $\widehat\prim$ be an arbitrary but fixed spread
of $\prim$, i.e.\ a choice of a translate of every circuit in $\prim$
at every vertex of $\gf$.

Start with the set $\circuits_0=\{\cdot\}$ consisting only of the
empty circuit; its growth series is $G_{\circuits_0}=1$. Let
$\widetilde{\circuits_1}$ be the saturation
$\langle\circuits_0\rangle$ of $\circuits_0$, and let $H(t,u)$ be its
spiky growth series, obtained via Theorem~\ref{thm:cogrowth}.

Let $\circuits_1$ be the set of circuits $c\in\widetilde{\circuits_1}$
such that, for all $p\in\widehat\prim$, the circuit $c$ does not
contain any prefix $p'$ of $p$ with $|p'|\ge(1-\eta)|p|$.

Since $f(t)$ counts the prime circuits $p$, we get $f(t^{1-\eta})$ as
the growth of ``forbidden'' prefixes $p'$. Solve $\zeta-1+(1-1/d)\zeta
f(\zeta^{\eta-1}/(d-1))=0$ for $\zeta$; then by
Lemma~\ref{lem:localtree} the spiky growth series of $\circuits_1$ is
bounded from below as
\[G_{\circuits_1}(t,u)\succsim H(\zeta t,u).\]

Consider next $\circuits_2$, the circuits obtained from $\circuits_1$ by
inserting at each vertex a non-negative number of prime circuits such
that the resulting circuit is reduced. Insertion of $0$ prime circuits
can be done in $0$ or $1$ ways, depending on whether the vertex is a
spike or not, and insertion of $i\ge1$ circuits can be done in at
least $(1-\frac2d)(1-\frac1d)^{i-1}f(t)^i$ ways, counting the insert's
length. Indeed, to guarantee that the resulting circuit is reduced, it
suffices to forbid one out of $d$ starting edges for all prime
circuits inserted, except for the last one, for which a starting and
an ending edge must be forbidden.

Summing over $i$ gives the generating functions
\[\frac{(d-2)f(t)}{d-(d-1)f(t)},\quad\frac{d-f(t)}{d-(d-1)f(t)}\]
counting possible insertions at a spike and non-spike vertex
respectively. The growth series of $\circuits_2$ is therefore
minorized as
\[G_{\circuits_2}(t,1)\succsim G_{\circuits_1}\left(t\frac{d-f(t)}{d-(d-1)f(t)},\frac{(d-2)f(t)}{d-f(t)}\right).\]
Finally let $\circuits_3=\langle\circuits_2\rangle$; a last application
of Theorem~\ref{thm:cogrowth} gives a lower bound for
$G_{\circuits_3}(t,1)$, which in turn, using
Proposition~\ref{prop:lowerbd}, gives a lower bound on $\|M\|$.

\begin{proof}[Proof of Theorem~\ref{thm:main}]
  The proof relies on the construction of cactus trees described in
  this section. In a cactus tree, i.e.\ a path constructed as above,
  mark the edges with the alphabet $\{1,2,3\}$ according to the stage
  at which that edge appeared in the cactus tree.
  
  The arguments given above show that the growth series
  $G_{\circuits_3}(t,1)$ undercounts marked cactus trees; indeed
  different choices of initial tree (in $\circuits_1$), prime circuits
  (in $\circuits_2$) or final spikes (in $\circuits_3$) yield
  different marked cactus trees. To show that $G_{\circuits_3}(t,1)$
  undercounts circuits, it therefore suffices to show that the
  markings on cactus trees are uniquely determined; i.e.\ that two
  distinct marked cactus trees remain distinct after the marks are
  erased.
  
  Consider a cactus tree. After removal of all spikes, it gives rise
  to a unique reduced path; in other words, the order in which the
  spikes are removed does not change the resulting reduced path.
  
  Now, in this reduced circuit, locate the first subword in
  $\widehat\prim$, and pluck it; and repeat till no such subword can
  be removed. These subwords are necessarily the prime circuits
  that were inserted in constructing $\circuits_2$; the only other
  possibility would be that some circuit $p=p'p''\in\circuits_1$ is
  such that $p'qp''$ contains a prime circuit $r$ at a position
  before $q$. This $r$ would then be either a subword of $p'$, which
  is forbidden by our construction of $\circuits_1$, or a subword of
  $p'q$ containing a part of $q$; by the small cancellation condition,
  a large part of $r$ would subsist in $p'$, and this is also
  forbidden by our construction.
\end{proof}

\section{Computations: Surface groups}
Here we make explicit the arguments in the previous section. Even
though our main motivation is to obtain lower bounds for the spectral
radius of surface groups, all computations are performed on $X_{d,m}$,
the $1$-skeleton of a tessellation of the hyperbolic plane by
$m$-gons, $d$ per vertex, introduced above.

Fix a base vertex $v$ of $X_{d,m}$. Consider as prime circuits the
$2d$ distinct $m$-gons touching $v$, in both orientations. Then
$\widehat\prim$ is the set of $m$-gons of $X_{d,m}$, taken with all
possible base points and orientations.

This choice amounts to taking $f(t)=2dt^m$. Since two $m$-gons touch
in at most one edge, we have a $\SC(1/m)$ small cancellation
condition, i.e.\ we take $\eta=1/m$.

By Equation~\eqref{eq:localsg}, we find $\zeta$ by solving
$\zeta-1+2(d\zeta-1)/(d^m\zeta^m-1)=0$ for $\zeta$'s largest
positive root.  Numerically, this root is close to $1$, and slightly
smaller; i.e.\ $\zeta=1-\mathcal O(e^{-d})$. For $d=8$ and $d=12$, we
obtain respectively
\[\zeta(8)=1-0.63\cdot10^{-5},\quad\zeta(12)=1-0.29\cdot10^{-10}.\]

We now follow the steps of Theorem~\ref{thm:main} for $d=m=8$,
corresponding to the surface group of genus $2$. We obtain
\begin{align*}
  h(t)&=\frac{2\cdot7}{6+8\sqrt{1-4\cdot7t^2}},\\
  H(t,u)&=\frac{1-(1-u)^2t^2}{1+(1-u)(7+u)t^2}h\left(\frac{t}{1+(1-u)(7+u)t^2}\right),\\
  g_1(t,u)&=H(t\zeta(8),u),\\
  g_2(t)&=g_1\left(t\frac{d-f(t)}{d-(d-1)f(t)},\frac{(d-2)f(t)}{d-f(t)}\right)\\
  &\cong\frac{14(1-t^2)(1-14t^8)}{(1+7t^2-14t^8-2t^{10})\left(6+8\sqrt{1-28\left(\frac{\zeta(8)t(1-2t^8)}{1+7t^2-14t^8-2t^{10}}\right)^2}\right)},\\
  g_3(t)&=h(t)g_2\left(\frac{1-\sqrt{1-28t^2}}{14t}\right).
\end{align*}
In particular, $g_1$ and $g_2$ are algebraic functions of degree $2$,
and $g_3$ is algebraic of degree $4$. Clearly $g_3$'s radius of
convergence is given by the $g_2(\cdot)$ term, since $g_2$ counts more
circuits than $h$; the radius of convergence of $g_2$ is given by
the vanishing of the surd expression
\[1-28\left(\frac{\zeta(8)t(1-2t^8)}{1+7t^2-14t^8-2t^{10}}\right)^2,\]
which is a polynomial equation of degree $20$, with least solution
$\alpha\approx0.357936$. Now the radius of convergence of $g_4$ is the
minimal $t$ such that $(1-\sqrt{1-28t^2})/14t=\alpha$, namely
$\rho=\alpha/(1+7\alpha^2)\approx0.188702$. We then have
$\|M_2\|\ge1/(8\rho)$, so
\[\|M_2\|\ge 0.662418.\]

Similar computations give
\[\|M_3\|\ge 0.552773.\]

The best upper bounds were obtained by Tatiana
Nagnibeda~\cite{nagnibeda:upperbd}, by applying Gabber's
lemma~\cite{gabber-g:superconcentrators} to a function on the edges of
$X_{4g,4g}$ defined by the edge-origin's cone type. She obtained
\[\|M_2\|\le 0.662816,\quad\|M_3\|\le 0.552792.\]

\subsection{Isoperimetric constants} Let $\gf$ be a connected
graph. For a subset $K\subset V$ of vertices denote by $\partial K =
\{e\in E:\,\alpha(e)\in K\text{ and }\omega(e)\not\in K\}$ the
\emph{boundary} of $K$. The number
\begin{equation}\label{iso}
  \iota(\gf) = \inf_{\text{finite, non-empty }K\subset V}\frac{\#\partial K}{\#K}
\end{equation}
is called the the \emph{(edge-)isoperimetric constant} of $\gf$. It is
connected to the spectral radius by the following
\begin{thm}[\cite{mohar:isoperimetric,mohar-w:spectra}]\label{thm:iotamu}
  Let $\gf$ be an infinite $d$-regular graph. Then one has
  \begin{equation}\label{eq:muiota}
    \frac{d^2-(d-2)\iota(\gf)}{d^2+\iota(\gf)}\leq\|M_\gf|\leq\sqrt{1-\iota(\gf)^2/d^2}.
  \end{equation}
\end{thm}
Note that, by~\eqref{eq:kesten} and \eqref{eq:muiota}, $\gf$ is
amenable if and only if $\iota(\gf) = 0$.

The isoperimetric constant $\iota(X_{d,m})$ has recently been computed
independently by Yusuke Higuchi and Tomoyuki
Shirai~\cite{higuchi-s:dfregular} and by Olle H\"aggstr\"om, Johan
Jonasson and Russell Lyons~\cite{haggstrom-j-l:isoperimetric}. They
obtained the values
\begin{equation}
  \iota(X_{d,m}) = (d-2)\sqrt{1 - \frac{4}{(d-2)(m-2)}}.
\end{equation}

In particular, they obtained $\iota(X_{8,8})=4\sqrt2$ which, together
with \eqref{eq:muiota} gives
$$0.431\approx\frac{16-6\sqrt2}{16+\sqrt2}\leq\|M_2\|\leq\frac{1}{\sqrt{2}}\approx0.707.$$
Evidently, these bounds are much weaker than the ones
given in this paper. It may be worthwhile to improve the connection
between the isoperimetric constant and the spectral radius ---
probably Theorem~\ref{thm:iotamu} is not the last word in this topic.

\bigskip\noindent{\bf Acknowledgment.}  The author thanks Tullio
Ceccherini-Silberstein, Pierre de la Harpe, Russ Lyons, Tatiana
Nagnibeda, and Yuval Peres, who provided valuable comments and
encouraged him to write this note.

\bigskip\noindent
Laurent Bartholdi,\\
Department of Mathematics.\\
University of California,\\
94720 Berkeley, USA

E-mail: laurent@math.berkeley.edu
\bibliography{mrabbrev,people,math,grigorchuk,bartholdi}
\end{document}